\documentclass{elsart}
\usepackage{amssymb,natbib,amsmath,enumerate,latexsym,array}

\usepackage[dvips]{graphicx} 
\usepackage[all]{xy} 
\newcommand{\mon} {\mathrm{mon}} 

\newcommand{\bZ}{\mathbb{Z}}
\newcommand{\cC} {\mathcal{C}} 
\newcommand{\cI} {\mathcal{I}}   
\newcommand{\cL} {\mathcal{L}}  
\newcommand{\cLinst} {\mathcal{L}_{inst}}  
\newcommand{\cLcomp} {\mathcal{L}_{comp}}  
\newcommand{\cP} {\mathcal{P}}

\newcommand{\src}{\mathrm{src}} 
\newcommand{\tgt}{\mathrm{tgt}} 
\newcommand{\GAP} {{\sf GAP}} 
\newcommand{\ep}{\epsilon}

\newcommand{\bA} {\mathsf{A}}
\newcommand{\bB} {\mathsf{B}}
\newcommand{\sets} {\mathsf{Sets}}
\newcommand{\arr}{\mathrm{Arr}} 
\newcommand{\ob}{\mathrm{Ob}} 
\def\geq{\geqslant}

\begin{document} 
\newtheorem{example}{Example}[section] 
\newtheorem{Alg}[example]{Algorithm}  
\newtheorem{Rem}[example]{Remark}  
\newtheorem{Lem}[example]{Lemma}  
\newtheorem{Thm}[example]{Theorem}  
\newtheorem{Cor}[example]{Corollary}  
\newtheorem{Def}[example]{Definition}  
\newtheorem{Prop}[example]{Proposition}  

\newenvironment{proof}{\noindent {\bf Proof} }{\mbox{} \hfill $\Box$ 
\mbox{}}  
\newenvironment{out}{\noindent {\bf Outline Proof} }{\mbox{} \hfill $\Box$ 
\mbox{}}

\begin{frontmatter}

\title{Logged Rewriting for Monoids}

\author[Leicester]{Anne Heyworth\thanksref{label1}}
\ead{a.heyworth@mcs.le.ac.uk},
\author[Macquarie]{Michael Johnson}
\ead{mike@ics.mq.edu.au},
\thanks[label1]{EPSRC GR/R29604/01 `Kan - A Categorical
Approach to Computer Algebra'.}
\address[Leicester]{Department of Computer Science,
Leicester University,
LE1 7RH, U.K.}
\address[Macquarie]{Information and Communication Sciences,
Macquarie University, NSW 2109
Australia.}

\title{}

\begin{abstract} 
 A rewriting system is a set of equations over a given set of terms called
rules that characterize a system of computation and is a powerful general
method for providing decision procedures of equational theories, based
upon the principle of replacing subterms of an expression with other
terms.
 In particular, a string rewriting system is usually associated with a
monoid presentation. At the first level the problem is to decide which
combinations of the generators are equivalent under the given rules;
Knuth-Bendix completion of the string rewriting system is one of the most
successful mechanisms for solving this problem. At the second level, the
problem involves determining which combinations of rules are equivalent.
 Logged rewriting is a technique which not only transforms strings but 
records the transformation in terms of the original system of rules. The 
relations between combinatorial, homotopical and homological finiteness 
conditions for monoids prompt us to consider using computer-friendly 
rewriting systems to calculate homotopical and homological structure from 
monoid presentations.
\end{abstract}


\end{frontmatter}

\section{Introduction} 

The idea of making a note of which rules are used as they are applied is
quite a simple one and it would be easy to regard it as too trivial to
spend any time on. However, when we look into the algebraic structure of
the records themselves, things become a lot less trivial. The background
to this work includes the papers by Squier, Lafont, Prout\'e, Otto,
Cremanns, Anick, Kobayashi, Pride and others on finiteness conditions for
monoids. The `combinatorial' finiteness condition is that a monoid has a
finite complete presentation. This implies the homological finiteness
condition $FP_\infty$~\citep{Anick,Squier87,Kobayashi} and also the
`homotopical' finiteness condition FDT \citep{Squier94}. It is also known
that the weaker homological condition $FP_3$ is not sufficient for either
FDT or the existence of a complete rewriting system in the case of
monoids~\citep{Squier94}.
 
Our aim is to use enhanced rewriting procedures to explicitly provide:
\begin{enumerate}[i)]
\item 
A finite complete rewriting system (combinatorial specification). 
\item 
A finite set of homotopy generators for $P^{(2)}\Gamma$ (homotopical 
specification). 
\item 
A (small) finitely generated resolution (homological specification). 
\end{enumerate} 
 
Logged rewriting for group presentations~\citep{paper5} gives procedures
for representing consequences of the relations of the presentation as
elements of a pre-crossed module and algorithms for computing generators
of the modules of identities among relations. In the monoid case the
relations are given by pairs of terms and the structure of a crossed
module is not appropriate to represent consequences of relations.  It is
well known \citep{Stell,Street} that sesquicategories or 2-categories can
be used to model rewriting systems. It has been
proved~\citep{Pride99,Gilbert} that when the rewriting system comes from a
group presentation, the 2-category can be identified with the
crossed module of the presentation.
 
In the special case of groups, various results are known.  In particular
there are methods for calculating a set of generators for the kernel
$\Pi_2$ of the crossed module of `consequences', which is useful for
constructing resolutions and calculating (co)homology.  For the case where
the rewriting system does not present a group we detail the algebraic
structure of the analogue of $\Pi_2$; presenting an algorithm for
computing a set of generators for it; and provide justification that the
constructions we make give combinatorial, homotopical and (co)homological
detail in the same spirit as $\Pi_2$.

In the case of monoids, logged rewriting techniques have further
applications. Specifically, we have so far examined applications to the
analysis of coset and double coset systems and used logged rewriting to
provide an alternative to the Reidemeister-Schreier algorithm for finding
presentations of subgroups~\citep{DCosets,RSAlt}. Additionally, we show in 
Section 8 that logged rewriting techniques are easily generalised to Kan 
extensions where they provide a proof technique for a wide range of 
decidablility problems solvable by string rewriting~\citep{paper2}.

\section{Logged Rewriting Systems} \label{LRS} 
 
A \emph{monoid presentation} is given as a pair $\mathcal{P} = \mon
\langle X | R \rangle$ where $X$ is a set of generators and $R$ is a set
of pairs $(l,r)$ of elements of the free monoid $X^*$. The monoid
presented, $M$ is the quotient obtained by factoring $X^*$ by $=_R$, the
congruence generated by $R$. The quotient monoid morphism will be denoted
$\theta:X^* \to M$.
 
We will assume that the reader is familiar with standard string rewriting
as in~\citep{BoOt}. The notation we use follows the usual conventions.  
The set $R$ is a \emph{rewriting system} for the monoid $M$ and its
elements are referred to as {\em rules}. The \emph{reduction relation
generated by $R$} on the free monoid $X^*$ is denoted by $\to_R$, the
reflexive, transitive closure is denoted $\stackrel{*}{\to}_R$ and the
reflexive, symmetric, transitive closure $\stackrel{*}{\leftrightarrow}_R$
coincides with the congruence $=_R$.  For convenience we assume that $R$
is compatible with an admissible well-ordering $>$; i.e. for all pairs
$(l,r) \in R$, we have $l>r$. This ensures that the relation $\to_R$ is
Noetherian.

The main aim of this section is to formally define a \emph{logged rewrite
system} for $\cP$. Such a system must not only reduce any word in $X^*$ to
an irreducible word (unique if the rewriting system is complete) but must
also express the actual reduction as a consequence of the original monoid
relations. The reader who does not wish to get into the details at this 
stage may wish to think of a consequence of the monoid relations as a 
sequence of rewrites recorded as: {\tt [prefix, rule, direction of rule, 
suffix]} which must give a valid rewrite.

It is important to identify the algebraic framework for these
`consequences' in order to understand what we may do with them. Formally,
one represents consequences of group relations by elements of a crossed
module. Consequences of monoid relations cannot be represented in that
framework; essentially this is because the free monoid does not have
inverses. However, it is well known that general string rewriting systems
may be modelled by sesquicategories or
2-categories~\citep{Benson,Stell,Street}. Therefore to every monoid
presentation we shall associate a sesquicategory. Its 2-cells correspond
to possible sequences of rewrites and inverse rewrites between strings in
the free monoid, with respect to the given rewriting system. Formally:

\begin{Def}[Sesquicategory of Rewrites]\mbox{ }\\
The \emph{sesquicategory $SQ(\cP)$ of a monoid presentation $\cP$} 
consists of the following:
\begin{itemize}
\item 
a single 0-cell which is denoted $*$,
\item
a free monoid of 1-cells which are the elements of $X^*$,
\item
a collection of 2-cells which are sequences
$$ 
\alpha = u_1 \alpha_1^{\ep_1}  v_1 \cdot 
  \cdots \cdot u_n \alpha_n^{\ep_n} v_n 
$$ 
where $u_1, \ldots, u_n, v_1, \ldots, v_n \in X^*$,  
$\alpha_1, \ldots, \alpha_n \in R \cup \{1\}$ and  
$\ep_1, \ldots, \ep_n = \pm 1$  
such that $u_i  \,   \tgt(\alpha_i^{\ep_i})   \,      v_i  
        =  u_{i+1} \, \src(\alpha_{i+1}^{\ep_{i+1}}) \, v_{i+1}$  
for $i=1, \ldots, n\! -\! 1$. 
\item
left and right actions of the 1-cells upon 
the 2-cells (whiskering) i.e. for any rewrite $\alpha$ and any 
elements $u$ and $v$ of the free monoid we say that $u \alpha v$ is a 
rewrite and 
 $\src(u \alpha v) = u \src(\alpha) v$ and 
 $\tgt(u \alpha v) = u \tgt(\alpha) v$.  
\item
identity rewrites for each string $w\in X^*$, denoted
$1_w$ where $\src(1_w)=\tgt(1_w)=w$
with the property that
$$ u \cdot 1_w \cdot v = 1_{uwv}$$ for all $u, v$ in $X^*$.
\item
a partial (`vertical') composition of rewrites, defined so that  
$\alpha \cdot \beta$ is a rewrite with  
$\src( \alpha \cdot \beta ) = \src( \alpha )$ and  
$\tgt( \alpha \cdot \beta ) = \tgt( \beta )$  
whenever $\tgt( \alpha ) = \src( \beta )$. 
\end{itemize}
\end{Def} 
 
For the above definition it can be verified that the sesquicategory axioms 
hold with respect to vertical composition and the whiskering action. 

Further, we shall allow rewrites to be cancelled by the reverse 
application of the rewriting sequence.  
The formal inverse of any rewrite $\alpha$ is
denoted $\alpha^{-1}$ where $\src(\alpha^{-1})=tgt(\alpha)$ and
$\tgt(\alpha^{-1})=\src(\alpha)$ and we allow that
$$ 
\alpha \cdot \alpha^{-1} = 1_{\src(\alpha)} \text{ for all rewrites } 
\alpha. 
$$ 
This gives $\cdot$ a groupoid structure, so we may refer to the 
sesquigroupoid $SQ(\cP)$. 

In the case where we can apply the rule $\alpha$ to one substring of a
string and the rule $\beta$ to another substring which is completely
disjoint from the first, it is natural to regard the order in which the
rules are actually applied as immaterial. This interchangability of
non-overlapping rewrites is captured by the interchange law on the
sesquicategory, giving us a 2-category. We shall denote the set of 2-cells
in $SQ(\cP)$ by $C_2$.

\begin{Def}[2-category of Rewrites]\mbox{ }\\
The 2-category of rewrites $C_2(\cP)$ is obtained by factoring the 2-cells 
$C_2$ of $SQ(\cP)$ by the interchange law:
$$ 
I= \{ (\alpha \, \src(\beta) \cdot \tgt(\alpha) \, \beta, 
\src(\alpha) \, \beta \cdot \alpha \, \tgt(\beta)) : \alpha, \beta \in C_2 
\}
$$
\end{Def} 

Specifically, the set of pairs of $I$ generates a relation on $C_2$
$$ \{ (\gamma \cdot u \alpha_1 v \cdot \delta,  
          \gamma \cdot u \alpha_2 v \cdot \delta : (\alpha_1, 
\alpha_2) \in I, u,v \in X^*, \gamma, delta \in C_2 \}
$$ 
and the reflexive, symmetric, transitive closure of this is $=_I$, which 
preserves both vertical composition and whiskering.
Congruence classes are formally denoted with square brackets so
$[\alpha]_I$ denotes the class of $C_2$ under $=_I$ that contains $\alpha$.
Whiskering and vertical composition are preserved and so may be applied to 
the congruence classes:
$u[\alpha]_Iv=[u \alpha v]_I$ for all $u,v \in X^*$ and $[\alpha]_I \cdot
[\beta]_I = [\alpha \cdot \beta]_I$.
A horizontal composition of the congruence classes may also be defined: 
$[\alpha]_I \circ [\beta]_I = [\alpha \, \src(\beta) \cdot \tgt(\alpha) 
\, \beta]_I$.

In the case of term rewriting one does not always wish to factor out by
the interchange law as it destroys the notion of length (number of steps)
of a rewrite. In the case of string rewriting we do not have to worry
about notions of length of derivation, thus we use the 2-category.
However, it should be noted that whilst rewrites may be represented
uniquely in the sesquicategory, the word problem for the 2-category is
generally unsolvable (generalisation of a crossed module). Like many, for
convenience, we abuse notation a little, representing rewrites that
should strictly be written as classes $[\alpha]_I$ by non-unique
representatives in the sesquicategory $\alpha$. So a pair of rewrites
$\alpha, \beta \in C_2$, are equivalent if and only if $[\alpha]_I =
[\beta]_I$.

In the context of groups, the sesquicategory associated to a monoid
presentation is well known. Pride proved that if the monoid presentation
involved is obtained from a group presentation then the associated
2-category is isomorpic (as a crossed module) to the free crossed module
associated to the group presentation~\citep{Gilbert,Pride99}. Logged 
rewriting for groups was established by using the crossed module structure 
for the logs.
We now formally define logged rewriting using the 2-category associated 
with a monoid presentation.

\begin{Def}[Logged Rewriting System]\mbox{ }\\
A \emph{logged rewriting system} for a presentation $\mathcal{P}$ of a 
monoid $M$ is a collection of 2-cells (rewrites)   
$$ 
\cL = \{\alpha_1, \ldots, \alpha_n\} 
$$ 
of the associated 2-category $C_2(\cP)$  so that the 
\emph{underlying rewriting system} 
$$ 
R_\cL = \{ (\src(\alpha_1), \tgt(\alpha_1) ), \ldots,    
        (\src(\alpha_n), \tgt(\alpha_n) )         \} 
$$ 
is a rewriting system for $M$. 
\end{Def}   

A rewriting system $R$ on a monoid $M$ generates a {\em reduction 
relation} $$\to_R = \{ (ulv,urv) : (l,r) \in R, u,v \in M \}. $$
The reflexive, transitive closure of this relation is denoted 
$\stackrel{*}{\to}_R$, and the reflexive, symmetric, transitive closure is 
denoted $\stackrel{*}{\leftrightarrow}_R$ and coincides with the 
congruence generated by $R$ on $M$, denoted $=_R$.
The logged reduction of a string by a rule $\alpha$ is written as: 
$u \alpha v: u \src(\alpha) v  \to u \tgt(\alpha) v$ and the rewrite
recorded is $u \alpha v$.

If the elements $(l,r)$ of a rewriting system on a free monoid $X^*$ are
ordered such that $l>r$ with respect to some well-ordering on $X^*$, then
the resulting reduction system is {\em Noetherian}; i.e. an irreducible
element is reached after finitely many reductions. A reduction system is
{\em confluent} if for any string $w$ there exists a unique irreducible
string $\bar{w}$ such that $w \stackrel{*}{\to_R} \bar{w}$. A rewrite
system is said to be {\em complete} if the corresponding reduction
relation is both Noetherian and confluent. This is a desirable property,
since any pair of strings $w_1, w_2$ can be reduced in a finite number of
steps to their irreducible forms $\bar{w_1}$ and $\bar{w_2}$ which will be
equal if and only if $w_1 =_R w_2$; i.e. the {\em word problem} is
decidable.

\section{Logged Completion} 

The Knuth-Bendix algorithm attempts to convert an arbitrary rewriting
system into a complete one by adding rules compatible with the ordering to
the system to try to force confluence. The key concept here is that of
{\em critical pairs} which are pairs of reductions which can be applied to
the same string to obtain two different results. The important critical
pairs are associated with the \emph{overlaps} of a the rules in the
rewriting system. When considering normal critical pairs we only care
about the sources and targets of the rewrites and the relevant information
identifying the overlap. When we are dealing with a logged rewriting
system it is necessary to think of the sequences of rules giving the
instructions permitting both of the rewrites and to include these logs as
part of the critical pair information.

\begin{Def}[Logged Critical Pairs]\mbox{ }\\
An overlap occurs between the logged rewrites 
$\alpha_1 :l_1 \to r_1$ and 
$\alpha_2: l_2 \to r_2$ of $\cL$ whenever one of the following is true: 
\vspace{-1ex}
\begin{center}
\begin{tabular}{llll} 
    i) \quad $u_1 l_1 v_1 = l_2$,
&  ii) \quad $u_1 l_1 = l_2 v_2$,
& iii) \quad $l_1 v_1 = u_2 l_2$,
&  iv) \quad $l_1 = u_2 l_2 v_2$. 
\end{tabular}
\end{center} 
for some $u_1, u_2, v_1, v_2 \in X^*$. 
The {\em logged critical pair} resulting from the overlap is a whiskered 
pair $(u_1 \alpha_1 v_1, u_2 \alpha_2 v_2)$ for the appropriate $u_1, 
u_2, v_1, v_2 \in X^*$.
\end{Def}
 
Given a monoid presentation $\mathcal{P}= mon \langle X | R \rangle$ we
can associate to it an \emph{initial logged rewriting system} $\cL_{init}$
which consists of one 2-cell $\alpha$ for each rule $(l,r)$ of $R$ with
$\src(\alpha)=l$ and $\tgt(\alpha)=r$. These 2-cells are the generators of
sesquigroupoid associated to the presentation.
 
If the initial logged rewriting system is not complete then we can attempt
to transform it into a complete logged rewriting system, by adding 2-cells
which will make the underlying rewriting system complete, in a version of
the Knuth-Bendix algorithm which records information that is usually
discarded. Clearly, this recorded completion terminates exactly when the
usual completion procedure would terminate.

\begin{Alg}[Logged Knuth-Bendix Procedure]\mbox{ }\\
\vspace{-4ex}
\begin{small}
\begin{center}
\begin{tabular}{p{13.5cm}}
\hline
\begin{enumerate}[{LKB}1:]  
\item
(Input) 
Let $\cP$ be a presentation of a monoid with generators $X$ and relations
$(l_1,r_1), \ldots, (l_n,r_n)$ where $l_1>r_1, \ldots, l_n>r_n$ for some 
well-ordering on the free monoid $X^*$.  Define $\cL_{init}$ to be the 
set of 2-cells or logged rules $\{\alpha_1, \ldots, \alpha_n\}$, where
$\src(\alpha_i)=l_i$ and $\tgt(\alpha_i)=r_i$ for $i=1, \ldots, n$.
\item 
(Initialise) 
Set $\cL_{all}=\cL_{init}$; $\cL_{new}=\cL_{init}$; and let $C$ be the 
empty list. 
\item
(Search for Overlaps and Record Critical Pairs) 
Whenever an overlap occurs between the rewrites $\alpha_a \in \cL_{all}$ 
and $\alpha_n \in \cL_{new}$, record the associated critical pair by 
adding the element $(u_a, \alpha_a, v_a, u_n, \alpha_n, v_n)$ to the list 
$C$ (where $u_a, v_a, u_n$ or $v_n$ may be the identity element). 
\item
(Attempt to Resolve Critical Pairs) 
Set $\cL_{new}=\emptyset$. 
For every element of $C$ consider the pair $(u_a r_a v_a, u_n r_n v_n)$, 
reducing each string by $\cL_{all}$ to the irreducible strings $z_a$ and 
$z_n$ respectively. If $z_a=z_n$ then the critical pair is said to resolve 
and it can be removed from $C$. Otherwise we must add a new logged rule to 
the system. If $\beta_a$ and $\beta_n$ are the logs of the reductions 
to $z_a$ and $z_n$ then the new logged rule is   
$\gamma=\beta_a^{-1} \cdot u_a \alpha_a^{-1} v_a \cdot u_n \alpha_n 
v_n \cdot \beta_n$ 
if $z_a > z_n$ and  
$\gamma=\beta_n^{-1} \cdot u_n \alpha_n^{-1} v_n  \cdot u_a \alpha_a 
v_a \cdot \beta_a$ if $z_n>z_a$. 
Add $\gamma$ to $\cL_{new}$ and to $\cL_{all}$. 
\item
(Loop) 
If $\cL_{new}$ is non-empty then loop to {LKB$3$}. 
Otherwise the procedure terminates: all critical pairs of $\cL_{all}$ have 
been tested, and resolve. 
\item
(Output) 
Output $\cL_{all}$, a complete logged rewriting system for $\cP$.
\end{enumerate}
\mbox{ }\\ 
\hline
\end{tabular}
\end{center}
\end{small}
\end{Alg} 

The immediate application for logged rewriting systems is in the provision 
of witnesses for computation. An ordinary complete rewriting system can 
determine whether or not two strings $s_1$ and $s_2$ represent the same 
element of the monoid; a logged rewriting system produces a proof in terms 
of a sequence of specific applications of the original monoid relations 
which will transform $s_1$ into $s_2$.
This is a fairly shallow application, although variations on it are useful 
in more complex algorithms such as~\citep{paper5}.

\section{Endorewrites}

Deeper information about the presentation can be gained by studying the
interaction of the relations with each other, known in group theory as the
{\em identities among relations}. The identities themselves represent
rewrite sequences which start at a word, send it through various
transformations and return it to its original form. For the monoid case,
we decided to refer to such rewrites as endorewrites. In the case of
monoids, the structure is necessarily less simple than the kernel of a
crossed module map. Note that we will continue to identify rewrites
which should strictly be written as classes $[\alpha]_I$ by (non-unique)
representatives in the sesquicategory $\alpha$. So $\alpha, \beta \in EQ
\subseteq C_2$, are equal as rewrites if and only if $[\alpha]_I =
[\beta]_I$.

\begin{Def}[Endorewrites]\mbox{ }\\
A 2-cell $\alpha \in C_2$ is an \emph{endorewrite} on a string $w$ if
$\src(\alpha)=\tgt(\alpha)=w$. 
\end{Def} 

The set of all endorewrites is actually the equaliser object of the two
maps $\src, \tgt : C_2 \to X^*$ in the category of sets. We denote it
$EQ$.

\begin{Lem}[Endorewrite Structure]\mbox{ }\\
The set of all endorewrites $EQ$ is the disjoint union of the sets $EQ_w$ 
for $w \in X^*$ where
$$EQ_w = \{ \alpha : \src(\alpha)=\tgt(\alpha)=w \}.$$ 
Each $EQ_w$ is closed under vertical composition; and their union $EQ$ is 
additionally closed under horizontal composition and whiskering.
\end{Lem}

Vertical composition is defined only within subsets $EQ_w$.  Horizontal
composition is defined across subsets: if $\alpha \in EQ_w$ and 
$\alpha' \in EQ_{w'}$, then $\alpha \circ \alpha' \in EQ_{ww'}$.  
Whiskering means that for any substring $s$ of a string $w$ there is an 
injective mapping $EQ_s \to EQ_w$ defined by $\alpha \mapsto u \alpha v$ 
where $usv=w$.

\begin{Lem}[Conjugate Endorewrites]\mbox{ }\\
If $\theta(w_1)=\theta(w_2)$ then for every $\beta \in C_2$ such that 
$\src(\beta)=w_1$ and $\tgt(\beta)=w_2$ there exists a bijection 
$\Phi_\beta:EQ_{w_1} \to EQ_{w_2}$ defined by 
$\alpha \mapsto \beta^{-1} \cdot \alpha \cdot \beta$. 
\end{Lem}

Thus the elements of $EQ_{w'}$ are all conjugates of elements of $EQ_{w}$
so it becomes logical that we should only seek generators for 
endorewrites $EQ_w$ of one representative string $w$ for each monoid 
element $\theta(w)$. The next lemma helps to make this concrete.

\begin{Lem}[Partial Action of Rewrites on Endorewrites]\mbox{ }\\  
There is a partial function $EQ \times (C_2, \cdot) \to EQ$,
defined by $\alpha^\beta = \beta^{-1} \cdot \alpha \cdot \beta$ 
for $\alpha \in EQ$ and $\beta \in C_2$ such that 
$\src(\alpha)=\tgt(\alpha)=\src(\beta)$. 
This satisfies the following properties:
\begin{enumerate}[i)]
\item
$\alpha^{1_{\src \alpha}} = \alpha$ for all $\alpha \in EQ$.
\item
$\alpha^{(\beta_1 \cdot \beta_2)} = (\alpha^{\beta_1})^{\beta_2}$ for all 
$\alpha, \beta_1, \beta_2 \in C_2$ such that 
$\src(\alpha)=\tgt(\alpha)=\src(\beta_1)$ and $\tgt(\beta_1)=\src(\beta_2)$.
\item
$u \alpha^\beta v = (u \alpha v)^{u \beta v}$ for all $u, v \in X^*$ 
whenever $\alpha^\beta$ is defined.  
\item
$(\alpha_1 \cdot \alpha_2)^\beta = \alpha_1^\beta \cdot \alpha_2^\beta$ 
for all $\alpha_1, \alpha_2 \in EQ$ and $\beta \in C_2$ such that 
$\tgt(\alpha_1 \circ \alpha_2)=\src(\beta)$. 
\end{enumerate}
\end{Lem}

The first two properties are the categorical equivalent of the properties
required for a partial monoid action, the second two show that the partial
action preserves the whiskering and vertical composition operations in
$EQ$. All the properties follow from the definitions of $(C_2, \cdot)$,
identity 2-cells $1_{\src(\alpha)}$ and the definition of $\alpha^\beta$.

Intuitively, $\alpha$ is like a circular walk: conjugating by $\beta$ just
means that we first walk down an additional path to the start of $\alpha$,
retracing our steps back along that path once the circular walk $\alpha$
is completed. Clearly the circular walk is not much more interesting for
having this initial path added to it and a guidebook that suggested all
conjugates of $\alpha$ were distinct jaunts would be absurd. Thus we
factor $EQ$ by this partial action and consider $\alpha$ to be 
equivalent to all its possible conjugates. Formally:

\begin{Lem}[Classes of Endorewrites]\mbox{ }\\
Let $C_2(\cP)$ be the 2-category of rewrites for a monoid presentation 
$\cP$ and let $EQ$ be the set of all endorewrites.
Then define 
$$J =  \{ (\alpha, \beta^{-1} \cdot \alpha \cdot \beta) : \alpha 
\in EQ, \beta \in C_2 \text{ and } \src \alpha = \src \beta \}.$$ 
Let $=_{I+J}$ be the smallest congruence on $E$ with respect to $\cdot$ 
and whiskering which contains both $J$ and the interchange law $I$.
Then the quotient $EQ^J = EQ/\! =_{I+J}$ is well-defined, preserving both 
vertical composition and whiskering.
\end{Lem}

To conclude this section we observe the following lemma.

\begin{Lem}[Structure of $EQ_w$]\mbox{ }\\ 
Let $C_2(\cP)$ be the 2-category of rewrites for a monoid presentation 
$\cP$ of a monoid $M$. 
Then $EQ_w^J$, the set of classes of endorewrites on any string $w \in X^*$,
is a $\mathbb{Z}M$-bimodule with respect to vertical composition and whiskering.
\end{Lem} 

\begin{proof}
Vertical composition of conjugacy classes of $EQ_w$ gives an abelian 
group structure: it is associative, with identity is $[1_w]_{I+J}$;
the inverse of $[\alpha]_{I+J}$ is $[\alpha^{-1}]_{I+J}$; 
and if $\alpha_1, \alpha_2 \in EQ_w$ then
$$[\alpha_1]_{I+J} \cdot [\alpha_2]_{I+J} =
  [\alpha_1 \cdot \alpha_2]_{I+J} =
  [\alpha_1 \cdot \alpha_2]_{I+J}^{\alpha_2^{-1}} =
  [\alpha_2 \cdot \alpha_1 \cdot \alpha_2 \cdot \alpha_2^{-1}]_{I+J} = 
  [\alpha_2 \cdot \alpha_1]_{I+J}.$$

The left and right whiskering actions of $X^*$ on $EQ_w^J$ restrict to 
well-defined left and right actions of $M$ since 
$u_1 \, \alpha \, v_1 =_{I+J} u_2 \, \alpha \, v_2$,
when $\theta(u_1)=\theta(u_2)$ and $\theta(v_1)=\theta(v_2)$ since:
$$ u_1 \, \alpha \, v_1 = 1_{u_1} \cdot \alpha \cdot 1_{v_1}
                        =_{I+J} 1_{u_2} \cdot \alpha \cdot 1_{v_2}
= u_2 \, \alpha \, v_2.$$
\end{proof}

\begin{Rem}
Note that horizontal composition is not abelian:
if $\alpha: w \to w$ and $\beta: z \to z$ then $\alpha \circ \beta: wz \to wz$ 
whilst $\beta \circ \alpha: zw \to zw$ 
and generally we cannot expect that $\theta(wz)=\theta(zw)$.
\end{Rem}

\section{Critical Pairs}

In this chapter we shall prove the intuitively reasonable idea that all
distinct circular routes come from examining the reconnection of
non-trivial diverging paths and thus provide a method for identifying all 
the interesting endorewrites of any completable rewriting system. 
Our main result requires that we first identify exactly what we mean by 
`a generating set of endorewrites'. 

A \emph{generating set} for $EQ$ must be a set of endorewrites $E$ such
that any other endorewrite of $EQ$ is equivalent under the interchange law
together with the conjugacy congruence $=_{I+J}$, to a product of
whiskered elements and inverse elements of $E$. Formally:

\begin{Def}[Generating Set for $EQ$]\mbox{ }\\ 
A \emph{generating set} for the endorewrites $EQ$ associated with a monoid 
presentation is a set $E \subseteq EQ$ such that
for any $\gamma \in EQ$ there exist $\alpha_1, \ldots, \alpha_n \in E$
such that
$$\gamma =_{I+J} u_1 \alpha_1 v_1 \cdot 
  \cdots \cdot u_n \alpha_n  v_n$$
for some $u_1, \ldots, u_n, v_1, \ldots, v_n \in X^*$ and
$\ep_1, \ldots, \ep_n \in \{-1,1\} $.
\end{Def}

Our main theorem claims that a set of generating endorewrites $E$, 
can be produced from the critical pairs which result from overlaps of 
the completed rewriting system. 
In order to prove the theorem we use digraph arguments, a digraph being 
associated with each endorewrite coming from a critical pair in the 
following way:


\begin{Lem}[Digraphs associated with Endorewrites]\label{lem-yield}\mbox{ }\\
Given two strings $w$ and $z$, any pair of logged reductions $\alpha_1, 
\alpha_2 : w  \to z$ is represented by a labelled digraph
which is associated uniquely with an endorewrite. 
\end{Lem}

\begin{proof}

\begin{tabular}{@{}p{12cm}@{}c}
Let $\alpha_1, \alpha_2 : w \to z$ in $C_2$.  Then we
have a digraph $D(\alpha_1,\alpha_2)$, as shown, and the associated
endorewrite $\delta(\alpha_1, \alpha_2) = \alpha_1 \alpha_2^{-1} \in
EQ_w$ can be obtained by reading the labels anticlockwise from the edges,
beginning at the vertex which is greatest with respect to $>$ on $X^*$.
&
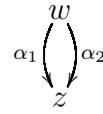
$$
\xymatrix{w \ar@/_/[d]_{\alpha_1} \ar@/^/[d]^{\alpha_2} \\ z}
$$
\end{tabular}
\end{proof}

\begin{Rem}[Resolved Critical Pairs Yield Endorewrites]\mbox{ }\\
If $C$ is the set of all logged critical pairs and $EQ$ is the set of all 
endorewrites of a complete logged rewriting system, then there is a 
map $\delta:C \to EQ$ associating an endorewrite with each critical pair.
In detail, if $c=(\alpha_1,\alpha_2)$. is a logged critical pair then there 
is a string $w$ which may be rewritten in two ways -- $\alpha_1:w \to w_1$ 
and $\alpha_2:w \to w_2$, where $\alpha_1, \alpha_2 \in C_2$.  
Since the pair can be resolved, there exists a string $z$ so that  
$\beta_1:w_1 \to z$ and 
$\beta_2:w_2 \to z$, 
for some rewrite sequences $\beta_1, \beta_2$ in $C_2$.  
It is immediate that 
$\delta(c)=\alpha_1 \cdot \beta_1 \cdot \beta_2^{-1} \cdot 
\alpha_2^{-1}$ is an endorewrite on $w$.
\end{Rem}

We now observe that endorewrites resulting from critical pairs are 
trivial when the critical pair involves disjoint rules or a 
conjugate of the endorewrite resulting from the reduction of the minimal 
string on which the same overlap occurs.


\begin{Lem}[Overlaps and Endorewrites]\label{lem-overlap}\mbox{ }\\
If  $\alpha_1:l_1 \to r_1$  and  $\alpha_2: l_2 \to r_2$ are rules of a 
complete logged rewriting system $\cL$, such that they may both be applied to 
a string $w$  then:
\begin{enumerate}[i)]
\item
if the rules overlap on $w$ then the endorewrite of the critical pair is
equivalent to a whiskering of the endorewrite given by a resolution of the 
same pair of rules applied to the minimal string on which the same 
overlap occurs.
\item 
if the rules do not overlap on $w$ then resolution of the critical pair
yields the trivial identity,
\end{enumerate}
\end{Lem}

\begin{proof}

\noindent
\begin{tabular}{@{}p{8.2cm}@{}l}
In case (i) the rules overlap on $w$ 
so there exist  
$u_1, v_1, v_2, x, y, z \in X^*$ 
such that  $w = xyz$  and either  
$y = u_1l_1v_1 = l_2$  or  
$y = u_1l_1    = l_2v_2$. 
In either case we can write  $y = u_1l_1v_1 = l_2v_2$ 
and the logged reductions of  $y$  are 
$u_1 \alpha_1 v_1: y \to u_1 r_1 v_1$ and  
$\alpha_2 v_2 : y \to r_2 v_2$.
By completeness there are logged reductions 
$\beta_1 : u_1 r_1 v_1 \to t$  and
$\beta_2 : r_2 v_2 \to t$ 
such that 
$\gamma = 
u_1 \alpha_1 v_1 \cdot \beta_1 \cdot \beta_2^{-1} \cdot \alpha_2^{-1} u_2$  
is an endorewrite.
The critical pair of reductions on  $w$  are
$x u_1 \alpha_1 v_1 z : w \to x u_1 r_1 v_1 z$ and
$x \alpha_2 u_2 z     : w \to x r_2 v_2 z$.
This pair can be resolved by
$x \beta_1 z : x u_1 r_1 v_1 z \to x t z$ and 
$x \beta_2 z : x r_2 v_2 z \to x t z$.
The endorewrite associated to it is 
 $x \gamma z$.
&
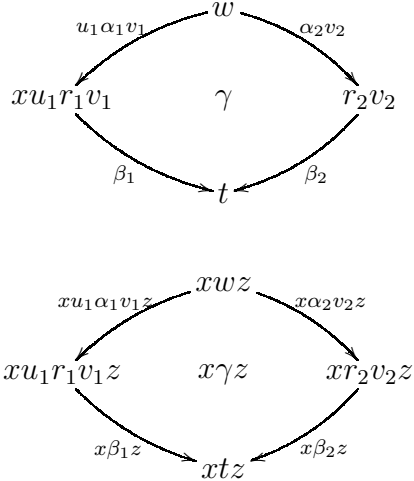
$$
\xymatrix{
   & w \ar@/_/[dl]_{u_1 \alpha_1 v_1} 
       \ar@/^/[dr]^{\alpha_2 v_2}&\\ 
 x u_1 r_1 v_1  \ar@/_/[dr]_{\beta_1} &  
\gamma 
&  r_2 v_2   \ar@/^/[dl]^{\beta_2}\\
  & t  &\\
   & xwz \ar@/_/[dl]_{x u_1 \alpha_1 v_1 z} 
         \ar@/^/[dr]^{x \alpha_2 v_2 z}&\\ 
 x u_1 r_1 v_1 z \ar@/_/[dr]_{x \beta_1 z} &  
  x \gamma z
  & x r_2 v_2 z  \ar@/^/[dl]^{x \beta_2 z}\\
  & x t z &\\}
$$
\end{tabular}

\noindent
\begin{tabular}{@{}p{7.8cm}@{}l}
In case (ii) the rules do not overlap on $w$ so there exist  $x,y,z \in X^*$ 
such that  $w = x l_1 y l_2 z$  and the logged reductions 
shown in the digraph on the right apply.
This yields the endorewrite
$x \alpha_1 y l_2 z \cdot x r_1 y \alpha_2 z \cdot 
 x l_1 y \alpha_2 z \cdot x \alpha_1 y r_2 z$, 
which is equivalent under the interchange law to $1_w$.
&
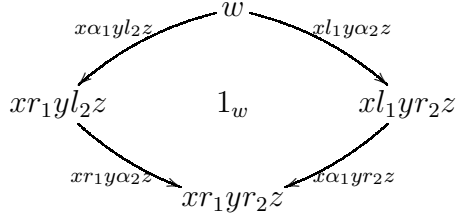
$$
\xymatrix{
   & w \ar@/_/[dl]_{x\alpha_1 y l_2 z} \ar@/^/[dr]^{x l_1 y \alpha_2 z}&\\ 
 x  r_1 y l_2 z  \ar@/_/[dr]_{x r_1 y \alpha_2 z} &  
1_w 
& x l_1 y r_2 z   \ar@/^/[dl]^{x \alpha_1 y r_2 z}\\
  & x r_1 y r_2 z &\\
 }
$$
\end{tabular}

\end{proof}


\begin{Lem}[Digraph of Reduction Sequences]\label{alg-graph}\mbox{ }\\
For any critical pair of logged reduction sequences, there exists a finite 
digraph which is the union of digraphs resulting from resolving critical 
pairs as in Lemma \ref{lem-overlap}.
\end{Lem}

\begin{proof}
Given two logged reduction sequences
$\alpha : w \to w_1 \to \cdots \to w_m \to z$ and 
$\alpha' : w  \to w_{m+1} \to \cdots \to w_n \to z$, 
we define a digraph $D$. 
The vertices $V(D)$ are the distinct words occurring in these sequences, and
there is an edge labelled $\alpha_i$ from $w_i$ to $w_j$ if $w_i \to w_j$
is a reduction step labelled by $\alpha_i$ in one of the two given 
reduction sequences.
The pair of reduction sequences $(\alpha, \alpha')$ yield the endorewrite 
$\gamma=\delta(\alpha, \alpha')$ in the 
way described in Lemma \ref{lem-overlap}.
We now add to the graph (if the graph is drawn, this looks like
subdivison into small confluence diagrams, the proof was originally
phrased in `diamonds').
Note that the vertices are ordered with respect to $>$ in $X^*$.

\noindent
\emph{
{\em \bf Algorithm 5.5 (Digraph Filling/Construction) }
\begin{small}
\begin{center}
\begin{tabular}{@{}p{13cm}}
\hline
\begin{enumerate}[{D}1:]
\item 
(Initialise)
Given $D$ as defined above, set $V$ to be the set of vertices in $D$ and set $i=1$. 
\item 
(Select a Vertex)
If $V$ is empty, go to step D7. Otherwise,
set $v_i$ to be the maximum vertex in $V$ and remove $v_i$ from $V$.
\item 
(Test and Resolve)
If the vertex is not the source of two distinct arrows in D then discard it and go back to step D2.
Otherwise, consider the corresponding two reductions 
$\beta_{i,1}: v_i \to v_{i,1}$ and $\beta_{i,2}: v_i \to v_{i,2}$ 
The critical pair $(\beta_{i,1},\beta_{i,2})$ 
can be resolved since $\cL$ is a complete rewrite
system so we have 
$\gamma_{i,1}: v_{i,1} \to z_i$ and 
$\gamma_{i,2}: v_{i,2} \to z_i$
\item 
(Create New Digraph)
Define $D_i$ to be the digraph 
$
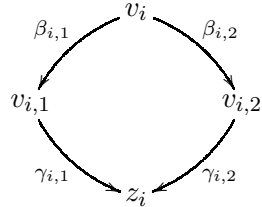
\xymatrix{
   & v_i \ar@/_/[dl]_{\beta_{i,1}} \ar@/^/[dr]^{\beta_{i,2}}&\\ 
 v_{i,1}  \ar@/_/[dr]_{\gamma_{i,1}} &  
& v_{i,2}    \ar@/^/[dl]^{\gamma_{i,2}}\\
  & z_i &\\
 }
$
\item
(Add to Digraph)
Add $D_i$ to $D$, identifying the vertices which have the same labels.
\item
(Loop)
Increment $i$ by 1 and go to step D2.
\item
(Terminate)
Output $D$.
\end{enumerate}\\
\hline
\end{tabular}
\end{center}
\end{small}
}

\vspace{1em}
We note firstly that $\cL$ is finite, so there are only 
finitely many rules which can be applied; 
secondly, any finite word can only be reduced in a finite number of ways; 
finally, the system is noetherian, so there are no 
infinite reduction sequences.
This means that the procedure will terminate, giving a finite digraph
$D$ which is the union of the digraphs $D_i$, which are all of the type
considered in Lemma \ref{lem-overlap}.
\end{proof}


\begin{Lem}[Digraph Compositions]\label{lem-graph}\mbox{ }\\
The product (at the base point) of the endorewrites associated (in the
sense of Lemma \ref{lem-yield}) with the sub-digraphs is equivalent under
the interchange law to the endorewrite associated with the original
digraph.
\end{Lem}
\begin{proof}

\begin{tabular}{@{}p{10.5cm}@{}l}
Consider the composition of 
digraphs of the type described, remembering that
each edge is associated uniquely to a particular log of the reduction.

The endorewrites associated to the two digraphs
are $\alpha_1 \cdot \gamma_1^{-1} \cdot \beta_1^{-1}$ 
and $\gamma_1 \cdot \alpha_2 \cdot \beta_2^{-1}$.
Composing them from the base point $w$ gives us
$\alpha_1 \cdot \gamma_1^{-1} \cdot \beta_1^{-1} \cdot 
\beta_1 \cdot (\gamma_1 \cdot \alpha_2 \cdot \beta_2^{-1}) \cdot 
                                                       \beta_1^{-1}$  
which is equivalent in the sesquigroupoid to
$\alpha_1 \cdot \alpha_2 \cdot \beta_2^{-1} \cdot \beta_1^{-1}$, the 
endorewrite given by taking the boundary of the composite. 
The fact that the order of the digraph endorewrites is not important 
corresponds with the fact that $EQ$ is abelian. 
&
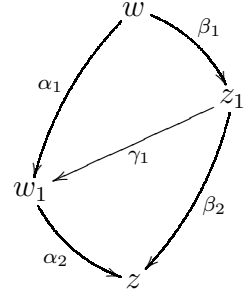
$$
\xymatrix{ & w \ar@/_/[ddl]_{\alpha_1} \ar@/^/[dr]^{\beta_1} &\\
           && z_1 \ar@/^/[ddl]^{\beta_2} \ar[dll]^{\gamma_1} \\
          w_1 \ar@/_/[dr]_{\alpha_2} && \\
         & z &}
$$
\end{tabular}

\end{proof}

Combining Lemma \ref{alg-graph} and Lemma \ref{lem-graph} with Lemma 
\ref{lem-overlap} we can deduce that any digraph can be identified with a  
product of whiskered endorewrites and inverse endorewrites of $E$. This 
allows us to prove the main theorem:

\begin{Thm}[Critical Pairs give a Set of Generators for 
$EQ$]\label{thm-main}\mbox{ }\\ 
Let $\cL_{init}$ be the initial logged rewriting system for a monoid 
presentation, and let $\cL_{comp}$ be a completion.
Let $C$ be the set of all logged critical pairs resulting from overlaps 
of the system $\cL_{comp} \cup \cL_{init}$. Then
$$E = \{ \delta(c) : c \in C \}$$
is a generating set of endorewrites.
\end{Thm} 

\begin{proof}
Let $\gamma$ be an endorewrite on some string $w$. 
Then consider the critical pair $(\gamma, 1_w)$.
Using Algorithm \ref{alg-graph} we can construct a digraph $D$ whose
associated endorewrite is $\gamma$ and whose sub-digraphs yield a product
of whiskered elements of $E$ and their inverses which is equivalent to 
$\gamma$ by Lemma \ref{lem-graph}.
\end{proof}

\section{Example}

This small example illustrates our methods for computing a complete set of
generators for the endorewrites of a monoid presentation from the overlaps
of a complete logged rewriting system.

Consider the monoid presentation $$mon\langle e,s \ | \ e^2=e, s^3=s, 
s^2e=e, es^2=e, sese=ese, eses=ese\rangle.$$ 
Using the short-lex ordering with $s>e$, labelling the relations
$\alpha_1, \ldots, \alpha_6$ we have the complete logged rewriting system 
consisting of the following six rules:

\begin{center}
\begin{tabular}{lll}
$\alpha_1 : e^2 \to e$, 
&
$\alpha_2 : s^3 \to  s$,
&
$\alpha_3 : s^2e \to  e$,
\\
$\alpha_4 : e^2s \to  e$,
&
$\alpha_5 : sese \to  ese$,
&
$\alpha_6 : eses \to  ese$.
\\
\end{tabular}
\end{center}

Consider the overlap of $\alpha_2$ and $\alpha_3$ on the string $w=s^3 e$.
Reducing it by $\alpha_2 e$ we get $se$ which is irreducible. 
Alternately, we can reduce $w$ by $s \alpha_3$ and similarly get $se$. 
Thus we have an endorewrite of $se$, i.e. 
$\alpha_2 e \cdot s \alpha_3^{-1}$.
Continuing in this way, considering all the overlaps of the logged system 
the following twenty six endorewrites can be computed:

Endorewrites of $e$: \;
$\alpha_2 se \cdot s^2 \alpha_3^{-1}$, \;
$\alpha_1 s^2 \cdot \alpha_4 \cdot \alpha_1^{-1} \cdot e \alpha_4^{-1}$, \;
$\alpha_1 e \cdot e \alpha_1^{-1}$, \;
$\alpha_3 e \cdot \alpha_1 \cdot \alpha_3^{-1} \cdot s^2 
\alpha_1^{-1}$, \;
$\alpha_3 s^2 \cdot \alpha_4 \cdot \alpha_3^{-1} \cdot s^2 
\alpha_4^{-1}$, \;
$\alpha_4 s^2 \cdot es \alpha_2^{-1}$ and
$\alpha_4 e \cdot e \alpha_3^{-1}$.

Endorewrites of $s$: \;
$\alpha_2 s^2 \cdot s^2 \alpha_2^{-1}$. 

Endorewrites of $s^2$: \;
$\alpha_2 s \cdot s \alpha_2^{-1}$.

Endorewrites of $es$: \;
$\alpha_4 s \cdot e \alpha_2^{-1}$.

Endorewrites of $se$: \;
$\alpha_2 e \cdot s \alpha_3^{-1}$.

Endorewrites of $ese$: \;
$\alpha_1 ses \cdot \alpha_6 \cdot \alpha_1^{-1} se \cdot e \alpha_6^{-1}$, \;
$s \alpha_5 \cdot \alpha_5 \cdot \alpha_3^{-1}se$, \;
$\alpha_3 ses \cdot s^2 \alpha_6^{-1}$, \;
$\alpha_4 se \cdot es \alpha_3^{-1}$, \;
$\alpha_5 e \cdot ses \alpha_1^{-1}$, \;
$\alpha_5 s^2 \cdot es \alpha_4 \cdot \alpha_5^{-1} \cdot ses \alpha_4^{-1}$, \;
$\alpha_5 s \cdot \alpha_6 \cdot \alpha_5^{-1} \cdot s \alpha_6^{-1}$, \;
$\alpha_5 ses \cdot \alpha_6 es \cdot es \alpha_1 s \cdot \alpha_6 
 \cdot \alpha_5^{-1} \cdot s \alpha_1^{-1} se \cdot se \alpha_5^{-1} 
 \cdot ses \alpha_6^{-1}$, \;
$\alpha_6 se \cdot \alpha_6 e \cdot ese \alpha_3^{-1}$, \;
$\alpha_6 s^2 \cdot es \alpha_4 \cdot \alpha_6^{-1} \cdot ese \alpha_2^{-1}$, \;
$\alpha_6 s \cdot \alpha_6 \cdot es \alpha_4^{-1}$, \;
$\alpha_6 e \cdot es \alpha_1 \cdot \alpha_1^{-1}se \cdot e \alpha_5^{-1}$, \;
$\alpha_6 ese \cdot ese \alpha_5$, \;
$\alpha_6 es \cdot es\alpha_1 s \cdot \alpha_6 \cdot es 
\alpha_1^{-1} \cdot \alpha_6^{-1} e \cdot es \alpha_6^{-1}$ and
$\alpha_5 se \cdot e \alpha_5 \cdot \alpha_1 se \cdot \alpha_5^{-1} 
\cdot s \alpha_1^{-1} se \cdot se \alpha_5^{-1}$.

These endorewrites generate all possible endorewrites of the system, but 
we note that generating sets obtained in this way are unlikely to be 
minimal generating sets.
For example, in this case there is a relation between the three 
endorewrites $\alpha_2 s \cdot s \alpha_2^{-1}$, \, $\alpha_2 e \cdot s 
\alpha_3^{-1}$, and $\alpha_2 se \cdot s^2 \alpha_3^{-1}$,
in that the third can be obtained from the first two in the following way:
$$ (\alpha_2 s \cdot s \alpha_2^{-1})e \cdot s(\alpha_2 e \cdot s 
\alpha_3^{-1}) = \alpha_2 se \cdot s^2 \alpha_3^{-1}.$$
Unfortunately, the fact that this problem generalises the word problem 
for crossed modules means that reducing the generating set can be rather 
ad-hoc since there are no normal forms for the 2-cells.

\section{Homotopical and Homological Interpretations}

We promised, in the introduction, that our results would enable
homotopical and homological specifications of the monoid. It is well
known that the existence of a finite complete rewriting system for a
monoid presentation implies the homological finiteness conditions
FP${}_3$~\citep{Squier87} and the stronger condition
FP${}_\infty$~\citep{Anick,Kobayashi} as well as the homotopical
condition of having finite derivation type
(FDT)~\citep{Cremanns,Squier94}. The addition made by this paper, in
considering logged rewriting systems, is that our algorithms enable the
specification of the structures which the properties are based upon.

In the homotopical case, it is immediate to observe that the set $E$ of
generating endorewrites suffices as a set of homotopy generators in the
sense of~\citep{Cremanns}. In detail: if $\alpha$ is any cycle of the
graph whose objects are all strings and whose invertible edges are all
rewrites, then $\alpha$ corresponds to the digraph of an endorewrite
and it turns out that the product of the subdigraphs is homotopically
equivalent to $\alpha$ for the same reasons as the associated
endorewrite is equivalent to the composite of the endorewrites of the
subdigraphs.

In terms of homology, the specification of $E$, similar to the analogous 
case of $\Pi_2$ for groups, enables us to construct a resolution. 
Specifically, we have an exact sequence of free, finitely generated 
$\mathbb{Z} M$-modules:

$$
\xymatrix{ 
   C_2 \ar[r]^{\delta_2} 
 & C_1 \ar[r]^{\delta_1} 
 & C_0 \ar[r]^{\delta_0}
 & \mathbb{Z} \ar[r]
 & 0.
}
$$

Given our specification of a finite set of homotopy generators, further
details of the resolution can be found in~\citep{Cremanns} in the proof
of the fact that FDT implies FP${}_3$.

For lower dimensional topology and cohomological dimensions for monoids,
Pride~\citep{Pride93,Pride95,Pride99} has developed geometric methods; 
using a calculus of
pictures, with spherical pictures representing the relations between the
relations, which may be identified with our endorewrites. His method for
determining a generating set differs significantly from ours; involving
picking an `obvious' set of pictures and then using picture operations to
prove that they generate all spherical pictures for the presentation.  
The key word here is `obvious' -- whether an obvious set of pictures can
be identified depends upon the shape of the presentation and its relation
to presentations for which generating sets of pictures are known. In the
case of groups substantial research means that many shapes of presentation
can be recognised, but in the case of monoids, presentations are less
recognisable.

Our generating set of endorewrites is determined algorithmically,
dependent on the successful completion of the presentation. The rewriting
method has the clear advantage of being able to be applied like brute
force in cases where the pictures are not obvious, or potentially in
complex examples where the pictures may be too complex to be identified by
eye. More interesting than comparing the two methods, however is to
consider using them in combination -- rewriting can provide an initial set
of pictures for unrecognisable monoid presentations and picture calculus
can then operate on the result to refine and reduce the set and present is
as something more ascetically pleasing and expressive than the strings of
letters representing whiskered 2-cells.

An alternative to looking at standard resolutions of a group by $\bZ
G$-modules as in \citep{Pride99} is to consider crossed resolutions. One
reason for interest in these is because their stronger invariance with
respect to the presentation makes them potentially more useful in the
classification of topological structures such as knots via crossed 
resolutions of their intertwining monoids.

Recall the group case: a {\em crossed complex (over groupoids)} is a 
sequence $C$
$$
\xymatrix{ 
 & \cdots \ar[r]^{\delta_{n+1}}
 & C_n \ar[r]^{\delta_n}
 & C_{n-1} \ar[r]^{\delta_{n-1}}
 & \cdots \ar[r]^{\delta_3}
 &  C_2 \ar[r]^{\delta_2} 
 & C_1 \ar[r]^{\delta_1} 
 & C_0 
}
$$
such that 
\begin{enumerate}[i)]
\item
$C_1$ is a groupoid with $C_0$ as its set of vertices and 
$\delta^1, \delta^0$ as its source and target maps.
\item
For $n \geq 2$, $C_n$ is a totally disconnected groupoid over $C_0$
and for $n \geq 3$, the groups at the vertices of $C_n$ are abelian.
\item
The groupoid $C_1$ operates on the right of each $C_n$ for $n \geq 2$ 
by an action denoted $(x,a) \mapsto x^a$.
\item
For $n \geq 2$, $\delta_n: C_n \to C_{n-1}$ is a morphism of 
groupoids over $C_0$ and $C_1$ acts on itself by conjugation.
\item
$\delta_{n}\delta_{n-1} = 0:C_n \to C_{n-2}$ for $n \geq 3$ and 
$\delta_2 \delta^0 
= \delta_2 \delta^1 : C_2 \to C_0$.
\item
If $c \in C_2$ then $\delta_2(c)$ operates trivially on $C_n$ for $n 
\geq 3$ and operates on $C_2$ by conjugation by $c$.
\end{enumerate}

A crossed complex $C$ is {\em free} if $C_1$ is a free groupoid (on some
graph $\Gamma_1$)  and $C_2$ is a free crossed $C_1$-module (for some
$\lambda: \Gamma_2 \to C_1$) and for $n \geq 3$, $C_n$ is a free
$\pi_1C$-module on some $\Gamma_n$ where $\pi_1 C$ is the fundamental 
groupoid of the crossed complex; i.e. the quotient of the groupoid $C_1$ 
by the normal, totally disconnected subgroupoid $\delta_2(C_2)$.

A crossed complex $C$ is {\em exact} if for $n \geq 2$ 
$$ Ker(\delta_n: C_n \to C_{n-1}) = Im(\delta_{n+1}:C_{n+1} \to C_n).$$

If $C$ is an free exact crossed complex and $G$ is a groupoid then $C$ 
together with an isomorphism $\pi_1C \to G$ (or, equivalently, C with a 
quotient morphism $C_1 \to G$ whose kernel is $\delta_2(C_2)$) is called a 
{\em crossed resolution of $G$}. It is a {\em free crossed resolution} 
if $C$ is also free.

In the case of monoids, we propose a similar structure.
Let $\cP = mon\langle X, R \rangle$ be a monoid presentation.
If we can find a complete rewriting system for $R$ then we can construct 
the following sequence:
$$
\xymatrix{ 
 & \cdots \ar[r]^{\delta_{n+1}}
 & C_n \ar[r]^{\delta_n}
 & C_{n-1} \ar[r]^{\delta_{n-1}}
 & \cdots \ar[r]^{\delta_3}
 &  C_2 \ar@<1ex>[r]^{\tgt} \ar@<-1ex>[r]^{\src} 
 & C_1 \ar[r]^{\delta_1} 
 & C_0 
}
$$

Define $C_0$ to be the monoid $M$ which is presented by $\cP$.
Define $C_1$ to be the free monoid $X^*$
and let $\delta_1:C_1 \to C_0$ be the quotient morphism. 
Then let $\src,tgt: C_2 \to C_1$ be the 
2-category of rewrites, but instead of 
a right action of $C_1$ we have a two-sided action; instead of a crossed 
module $\delta_2: C_2 \to C_1$ we have a 2-category $\src,\tgt:C_2 \to 
C_1$ and instead of $C_1$ being a groupoid, it is a category.
Then let $C_3$ be a family of free $\bZ M$-bimodules: 
its objects are the elements of $M$ and its arrows are of 
the form $\ep_1(m_1e_1n_1) + \ep_2(m_2e_2n_2) + \cdots \ep_k(m_ke_kn_k):m \to m$
when $m_1e_1n_1 \cdot m_2e_2n_2 \cdot \cdots \cdot m_ke_kn_k$ is an endorewrite
in $EQ_w$ for some $\theta(w)=m$.
For higher levels $n > 3$ we can define $C_n$ to be the free 
$\bZ M$-bimodule on a set of generators for $Ker(\delta_{n-1})$.

We find that $C$ is a crossed complex and we have maps $b_{i,j}: C_i \times C_j 
\to C_{i+j}$
-- whiskering in the case of $C_0$ operating on the left and right of 
$C_i$ for $i>0$. Then $C_1$ has 2 multiplications under the operations 
of $C_0$ which coincide only if $C_1$ is a monoid in the category of 
groupoids (interchange law).
There are no inverses in dimension 0, but inverses at all higher levels.
From the definitions we deduce exactness: $Ker(\delta_n) = 
Im(\delta_n+1)$.

This appears to be identifiable with the structure of a crossed
differential algebra, that is a crossed complex $C$ with a 
morphism $C \otimes C \to C$ which gives a monoid structure on $C$
(these are defined in detail in~\citep{Tonks}). 
We are still investigating how useful this enhanced style of resolution
may be in the monoid case, so we won't pursue the details of the 
construction further in this paper.

\section{Generalised Logged String Rewriting}

In~\citep{paper2} it was shown that the familiar string rewriting methods
can be applied to problems of computing left Kan extensions over the
category of sets. Structures such as monoid and category presentations,
induced actions of groups and monoids, equivalence and conjugacy classes,
equalisers and pushouts all turn out to be special cases of left Kan
extensions over $\sets$ and thus string rewriting methods can be
applied to all these variations on the word problem.

Since string rewriting for Kan extensions can be achieved by embedding
in a monoid, it is unnecessary to go through the detail of the
sesquigroupoid whose 2-cells possess the structure for the logged
rules. However, since we don't need to embed in a monoid in order for
the string rewriting methods to work, we briefly outline the alternative
sesquigroupoid.

Let $(E, \ep)$ be the left Kan extension of the category action $X: \bA 
\to \sets$ along the functor $F : \bA \to \bB$.
We assume that the data for the Kan extension is given as a finite 
presentation $\cP$, consisting of generating graphs for $\bA$ and $\bB$, a 
set of relations for $\bB$ and the action of functors $F$ and $X$ being 
defined for every object and arrow of the generating graph of $\bA$. 
The 2-category $C_2$ associated with the presentation of the Kan 
extension has 0-cells 
 $(\bigsqcup_{A \in \ob \bA} XA) \sqcup \ob \bB$ and 1-cells 
 $\{ (s_x:x \to FA)\, | \, x \in XA, A \in \ob \bA \} \sqcup \arr \bB$.
The 2-cells are the rewrites and inverse rewrites, with vertical 
composition as before, but clearly, whiskering and horizontal 
compositions are partial operations dependent on whether paths can be 
composed.

In conjunction with~\citep{paper2}, this observation enables logged
rewriting techniques to be applied to a wide range of problems, including
category presentations, equivalence relations, induced actions, pushouts
and coset systems. In eaach case, interpretations and potential
applications of the endorewrites requires further investigation.

\section{Implementations and Further Applications}

Techniques of logged rewriting have been implemented by the first author
as \GAP \, functions which will eventually be submitted as a package.
Applications of logged rewriting were explored in \citep{paper5} where the
group version was implemented, providing a new algorithmic method for the
construction of crossed resolutions of groups; in \citep{RSAlt} where the
logged completion methods give an alternative to the Reidemeister-Schreier
method of computing a subgroup presentation; and in \citep{DCosets} we
show how endorewrites for double coset rewriting systems reveal
information about the subgroups.

Further work could pursue other potential applications, including in Petri
nets, concurrency and the analysis of knot quandles; as well as
generalising the techniques to Gr\"obner bases where the endorewrites can
be identified with syzygies.


{\small 
\begin{thebibliography}{} 
 
\bibitem[Anick(1986)]{Anick} D. J. Anick,
                {\em On the Homology of Associative Algebras,} 
                Transactions of the American Mathematical Society, 
                 vol.296 p641-59 1986. 
 
\bibitem[Baader and Nipkow(1998)]{TAT} F. Baader and T. Nipkow,  
                {\em Term Rewriting and All That,} 
                Cambridge University Press 1998. 

\bibitem[Benson(1975)]{Benson} D.B. Benson,
                {\em The Basic Algebraic Structures in Categories of 
                 Derivations,} 
                Information and Control, vol. 28, no. 1, 
                Academic Press 1975.
 
\bibitem[Book and Otto(1993)]{BoOt}  R. Book and F. Otto,
                {\em String-Rewriting Systems,}  
                Springer-Verlag, New York 1993. 

\bibitem[Brown and Heyworth(2000)]{paper2} R. Brown and A. Heyworth,
\emph{Using Rewrite Systems to Compute Kan Extensions and Induced
  Actions of Categories},  J. Symbolic Computation 29 5-31, 2000.

\bibitem[Brown et al(2004)]{DCosets} R. Brown, N. Ghani, A. Heyworth and 
C.D. Wensley, {\em String Rewriting for Double Coset Systems,}
 submitted J. Symbolic Computation 2004.

\bibitem[Brown and Razak(1999)]{Brown} R. Brown and A Razak Salleh, 
                {\em Free Crossed Resolutions of Groups and Presentations of   
                Modules of Identities Among Relations,} 
               Journal of Mathematical Computation, LMS, vol.2 p28-61 
               1999.  
 
\bibitem[Buchberger and Winkler(1998)]{RISC} B. Buchberger and F. Winkler,  
                {\em Gr\"obner Bases and Applications,} 
                ``33 Years of Gr\"obner Bases'' RISC-Linz 2-4 Feb 1998, 
                Proc. London Math. Soc. vol.251. 
  
\bibitem[Cohen(1997)]{Cohen} D. Cohen,
                {\em String Rewriting and Homology of Monoids,} 
                 Math. Struct. Comput. Sci. 7 207-240 1997 
 
\bibitem[Cremanns(1995)]{Cremanns} R. Cremanns,
                {\em Finiteness Conditions for Rewriting Systems} 
                PhD thesis, Universit$\ddot{\text{a}}$t 
Gesamthochschule 
                Kassel, 1995.

\bibitem[Cremanns and Otto(1994)]{CrOt94} R. Cremanns and F. Otto,
                {\em Finite Derivation Type implies the Homological 
                     Finiteness Condition $FP_3$,} 
                Journal of Symbolic Computation, vol.18 p91-112 1994. 
 
\bibitem[Cremanns and Otto(1996)]{CrOt96} R. Cremanns and F. Otto,
                {\em Finite Derivation Type is equivalent to $FP_3$ for 
                Groups,} 
                Journal of Symbolic Computation, vol.22 p155-177 1996. 

\bibitem[GAP(1998)]{GAP}  
{{} The GAP~Group},  
\emph{GAP -- Groups, Algorithms, and Programming, Version 4}, Aachen, 
St~Andrews, 1998.
\verb+http://www.gap.dcs.st-and.ac.uk/~gap+. 

\bibitem[Ghani and Heyworth(2003)]{RSAlt} 
N. Ghani and A. Heyworth,
\emph{A Rewriting Alternative to the Reidemeister-Schreier Procedure},
RTA 2003.
 
\bibitem[Gilbert(1996)]{Gilbert} N. D. Gilbert,
                 {\em Monoid Presentations and Associated Groupoids,} 
                 Int. J. Algebra and Computation, 8 141-152 1998.
 
\bibitem[Groves(1997)]{Groves} J. R. J. Groves, 
                {\em An Algorithm for Computing Homology Groups,} 
                \emph{Journal of Algebra, vol.94 p331-361} 1997. 
      
\bibitem[Heyworth and Wensley(1999)]{paper5}
   {{}A. Heyworth and C.D. Wensley},
   {\em Logged Rewriting and Identities Among Relators},
   in Groups St Andrews 2001 in Oxford,
   eds. C.M. Campbell, E.F. Robertson, G.C. Smith,
   London Math. Soc. Lecture Note Ser. 304,
   C.U.P. p256-76 2003.
                                                           
\bibitem[Kobayashi(1990)]{Kobayashi} Y. Kobayashi,
         {\em Complete Rewriting Systems and Homology Of Monoid Algebras,}
         Journal of Pure and Applied Algebra, vol.65 1990 p263-275.

\bibitem[Lafont(1995)]{Lafont} Y. Lafont,
         {\em A New Finiteness Condition for Monoids Presented by
         Complete Rewriting Systems (after Craig C. Squier),}
         Journal of Pure and Applied Algebra, vol.98 1995 p229-244.

\bibitem[Lafont and Prout\'e(1991)]{LaPr} Y. Lafont and A. Prout\'e,
         {\em Church Rosser Property and Homology Of Monoids,}
         Mathematical Structures in Computer Science vol.1,
                  Cambridge University Press 1991 p297-326.

\bibitem[Mitchell(1972)]{Mitchell} B. Mitchell,
         {\em The 2-category of a Set of Relations,}
         in `Rings with several objects', Advances in Math. 8
         1972 p74-77.

\bibitem[Pride(1993)]{Pride93} S. J. Pride,
         {\em Low-dimensional Homotopy Theory for Monoids,}
         International Journal of Algebra and Computation,
         vol.5 1993 p631-649.

\bibitem[Pride(1995)]{Pride95} S. J. Pride,
         {\em Geometric Methods in Combinatorial Semigroup Theory,}
         in: Semigroups, Formal  Languages and Groups (J. Fountain ed), 
         Kluwer Academic Publishers, 1995 p215-32.

\bibitem[Pride(1999)]{Pride99} S. J. Pride,
         {\em Low-dimensional Homotopy Theory for Monoids II,}
         Glasgow Math Journal, 41 1999 p1-11.

\bibitem[Reinert(1995)]{Birgitsthesis} B. Reinert,
         {\em On Gr\"obner Bases in Monoid and Group Rings,}
         PhD Thesis, Universit\"at Kaiserslautern 1995. 
   
\bibitem[Reinert and Zecker(1998)]{MRC} B. Reinert and D. Zecker,
         {\em MRC - A System for Computing Gr\"obner Bases in
                Monoid and Group Rings,}
         Universit\"at Kaiserslautern Preprint 1998.

\bibitem[Sims(1994)]{Sims} C. C. Sims,
                {\em Computation with Finitely Presented Groups,} 
                Cambridge University Press 1994. 
 
 
\bibitem[Squier(1987)]{Squier87} C. C. Squier,
   {\em Word Problems and a Homological Finiteness Condition for Monoids,} 
        Journal of Pure and Applied Algebra, vol.49 p201-17 1987. 
 
\bibitem[Squier(1994)]{Squier94} C. C. Squier, F. Otto and Y Kobayashi,
                {\em A Finiteness Condition for Rewriting Systems,} 
                Theoretical Computer Science, vol.131 p271-94 1994. 
 
\bibitem[Stell(1994)]{Stell} J. G. Stell,
          {\em Modelling Term Rewriting Systems by Sesquicategories,} 
          Technical Report TR94-02, University of Keele 1994. 
 
\bibitem[Street(1992)]{Street} R. Street,
                {\em Categorical Structures,} 
                in `Handbook of Algebra', M. Hazelwinkel (ed) vol.1 1992. 

\bibitem[Tonks(1993)]{Tonks} A. Tonks,
          {\em Theory and Applications of Crossed Complexes,} 
          PhD thesis, University of Wales, Bangor 1993.

\end{thebibliography} 
}
\end{document}